\documentclass[a4paper,12pt, reqno]{amsart}

\usepackage{graphicx}
\usepackage{hyperref}
\usepackage{amsmath}
\usepackage{amssymb}
\usepackage{amsthm}
\usepackage{appendix}
\usepackage{verbatim}
\usepackage{framed}
\usepackage{mathtools}
\usepackage{lipsum}

\makeatletter
\renewcommand{\pod}[1]{\mathchoice
  {\allowbreak \if@display \mkern 18mu\else \mkern 8mu\fi (#1)}
  {\allowbreak \if@display \mkern 18mu\else \mkern 8mu\fi (#1)}
  {\mkern4mu(#1)}
  {\mkern4mu(#1)}
}

\allowdisplaybreaks[1]

\theoremstyle{plain}
\newtheorem{thm}{Theorem}

\newtheorem{cor}[thm]{Corollary}

\theoremstyle{definition}

\renewcommand{\th}{\textsuperscript{th} }

\newcommand{\cited}[1]{~\cite{{#1}}}
\newcommand{\longtitle}[1]{%
  \ifodd\value{page}%
  {
    \protect\parbox{0.97\linewidth}{#1}\hfill%
  }
  \else%
    \hfill\protect\parbox{0.97\linewidth}{#1}%
  \fi%
}

\begin{document}
\title{On the Density of Weak Polignac Numbers}
\author{Stijn S.C. Hanson}
\email{stijnhanson@gmail.com}
\maketitle

\begin{abstract}
  Let $k$ be an integer which is the difference between prime numbers infinitely often. It is known that there are infinitely many such $k$ and, in this paper, we give a new unconditional proof that these $k$ have positive density and improve on current bounds, assuming a strong hypothesis.
\end{abstract}

In 2013, Yitang Zhang proved the bounded gap conjecture\cited{Zhang} which asserted that there exists some natural number $k$ such that there exist infinitely many pairs of primes whose difference is precisely $k$. We call any $k$ that satisfies such a property a \emph{weak Polignac number}.

More precisely, we first ask what sorts of finite sets $\mathcal{H}$ we could construct such that their translates $n + \mathcal{H}$ could be prime infinitely often. That is to say, we want to create a set $\mathcal{H} = \{h_1, h_2, \ldots , h_k\}$ where, for all primes $p$, there is some natural number $m$ such that
\begin{equation*}
  h_i \not\equiv m \pmod{p}
\end{equation*}
for all $1 \leq i \leq k$. Any such $\mathcal{H}$ is called an \emph{admissible set}.

This theorem of Zhang's (which was later improved upon by James Maynard\cited{Maynard} and Polymath Project 8\cited{Polymath} does this by saying that the translates $n + \mathcal{H}$ of an admissible set $\mathcal{H}$ contain at least two primes infinitely often. In other words
\begin{thm}[Maynard--Polymath--Zhang, 2014]
  Suppose that \; \; \; \; \;$\mathcal{H} = \{h_1, h_2, \ldots , h_k\}$ is an admissible set and define the difference set
  \begin{equation*}
    \mathcal{D} = \{h_j - h_i: h_i < h_j\}.
  \end{equation*}

  Then, for $k = 50$, $\mathcal{D}$ contains a weak Polignac number. 
\end{thm}

If we assume stronger conditions such as the Elliott--Halberstam conjecture or the generalised Elliott--Halberstam conjecture then we get the above theorem with $k = 5$ and $k = 3$ respectively.

If we let $\mathcal{P}(x)$ be the number of weak Polignac numbers less than or equal to $x$ then Pintz has shown\cited{Pintz} that
\begin{equation*}
  \liminf_{x \rightarrow \infty} \frac{\mathcal{P}(x)}{x} > 0
\end{equation*}
and, in fact, he obtained a positive bound for the lower density. This paper gives an unconditional bound which, while not being as strong as Pintz', has a slightly simpler proof so that, if all we care about is that the lower density is positive, this method will suffice. We also give a bound which improves upon Pintz' under the assumption of the generalised Elliott--Halberstam conjecture and show that, if it is not optimal, then the methods used can improve upon the bound by no more than $1/36$.

The method used here is a packing problem: we construct an infinite sequence of admissible sets $\mathcal{H}^n$ such that the difference sets $\mathcal{D}^n$ are disjoint. Previous work in the field has looked for admissible sets of minimal diameter. These are very useful for lowering the bound on the Maynard--Polymath--Zhang theorem but are not very suited to this packing problem as they are rather irregular.

To that end we define a \emph{regular admissible set of size $k$} to be an admissible set of the form
\begin{equation*}
  \mathcal{H}_k ^n = \{0, nP(k), \ldots , (k - 1)nP(k)\}
\end{equation*}
where $P(k) = \prod_{p \leq k} p$ is the $k$\th primorial and note that their difference sets are
\begin{equation*}
  \mathcal{D}_k ^n = \{nP(k), \ldots , (k - 1)nP(k)\}.
\end{equation*}

Consider all regular admissible sets of size $k$ $(\mathcal{H}_k ^n)_{n \in \mathbb{N}}$. We want to find a subset $\mathcal{N} \subseteq \mathbb{N}$ such that the sequence $(\mathcal{D}^n)_{n \in \mathcal{N}}$ is pairwise disjoint and which pack into the interval $[1, x]$. 

We first note that there can be at most $\frac{x}{(k - 1)P(k)}$ different regular admissible sets which have difference sets which are entirely contained within $[1, x]$. Then, for $n < m$,  $\mathcal{D}_k ^n$ overlaps with all regular admissible sets $\mathcal{D}_k ^m$ such that
\begin{equation*}
  im = jn
\end{equation*}
whenever $1 \leq i < j \leq k - 1$. There are $\frac{1}{2}(k - 1)(k - 2)$ different choices of $i$ and $j$ so every regular admissible set included in our collection excludes no more than that many other potential admissible sets. Therefore
\begin{align*}
  \mathcal{P}(x) &\geq \left[\frac{x}{(k - 1)P(k)}\right] - \frac{(k - 1)(k - 2)}{2}\sum_{i = 1} ^x 1_{\mathcal{N}}(i) \nonumber \\
  &\geq \left[\frac{x}{(k - 1)P(k)}\right] - \frac{(k - 1)(k - 2)}{2}\mathcal{P}(x)
\end{align*}
which implies that.
\begin{equation*}
  \mathcal{P}(x) \geq \frac{2x}{(k - 1)((k - 1)(k - 2) + 2)P(k)} + O(1).
\end{equation*}
Dividing through by $x$ gives the following result.
\begin{thm}
  Suppose that all admissible sets of size $k$ contain a weak Polignac number in their difference sets. Then
  \begin{equation*}
    \liminf_{x \rightarrow \infty}\frac{\mathcal{P}(x)}{x} \geq \frac{2}{(k - 1)((k - 1)(k - 2) + 2)P(k)}.
  \end{equation*}
\end{thm}

Putting $k = 50$ into this formula shows that, unconditionally, the weak Polignac numbers have density greater than $\frac{1}{35,462,538,431,226,065,088,930} > 2.819 \times 10^{-23}$. Then we get a new derivation of a theorem of Pintz \cite{Pintz}
\begin{cor}
  The lower asymptotic density of the weak Polignac numbers is positive.
\end{cor}
Which, in turn, allows us to apply such wonderful theorems as that of Szemer\'{e}di's\cited{Szemeredi} which tells us that there are arbitrarily long arithmetic progressions within the prime numbers.

Furthermore, under Elliott--Halberstam ($k = 5$) the density is greater than $\frac{1}{840} > 0.00119$ and, under Generalised Elliott--Halberstam ($k = 3$), $\frac{1}{24} = 0.041\dot{6}$.

Actually, under the Generalised Elliott--Halberstam conjecture we can do a little better. Let
\begin{equation*}
  \mathcal{A} = \{6n \leq x - 2: n \in \mathbb{N}\}
\end{equation*}
and write $\mathcal{A} = \{a_1, a_2, \ldots , a_N\}$ where $a_i = 0$ whenever $3 \mid i$ and $a_i \neq 0$ otherwise. We can also insist that the non-zero $a_i$ are strictly decreasing and define
\begin{equation*}
  \mathcal{H}^n = \{0, 2n, 2n + a_n\},
\end{equation*}
where $3 \nmid n$ and $n \leq [x/6]$.

These sets are admissible (as all of their elements are congruent to $0$ modulo $2$ and either $0$ or $2$ modulo $3$). Also all the elements of the difference sets are contained within $[2, x]$ (as the largest element is $2 + a_1 \leq x$ by definition). Therefore, in this case,
\begin{equation*}
  \liminf_{x \rightarrow \infty} \frac{\mathcal{P}(x)}{x} \geq \frac{1}{6}.
\end{equation*}

We might also like to ask just how successful this stratagem can be. We know that the $k$-element admissible sets with the smallest difference sets are regular admissible sets with $k - 1$ elements so we can pack no more than
\begin{equation*}
  \left[\frac{x}{2(k - 1)}\right]
\end{equation*}
even difference sets into $[1, x]$. This gives a maximum packing density of $\frac{1}{2(k - 1)}$ assuming only the bounded gaps conjecture holds for $k$ and not for $k - 1$.

We can improve upon this in particular cases, however. Take $k = 3$ and note that the maximum cardinality for the difference set of an admissible set is $3$ and the minimum is $2$. Moreover, this minimum is achieved precisely when the admissible set is regular. We know that regular admissible sets can only be made up of multiples of $P(k) = 6$ and they will have $k - 1 = 2$ elements  so there are no more than
\begin{equation*}
  \left[\frac{x}{12}\right] \sim \frac{x}{12}
\end{equation*}
even regular admissible sets with disjoint difference sets contained in $[1, x]$. We can therefore have an extra
\begin{equation*}
  \left[\frac{\left[\frac{x}{2}\right] - 2\left[\frac{x}{12}\right]}{3}\right] \sim \frac{1}{9}
\end{equation*}
irregular admissible sets. This tells us that we can have no more than
\begin{equation*}
  \frac{x}{12} + \frac{x}{9} = \frac{7x}{36}
\end{equation*}
admissible sets whose difference sets are disjoint and contained in $[1, x]$ for large enough $x$. Therefore, only assuming the generalised Elliott--Halberstam conjecture, the maximum packing density is $7/36$.

\end{document}